\documentclass[1pt]{elsarticle}
\usepackage{geometry}
\usepackage{graphicx,subfigure}
\usepackage[colorlinks,linkcolor=red,citecolor=blue]{hyperref}
\usepackage{amsmath,amssymb,amsxtra,bm}
\usepackage{mathrsfs,amsfonts,amsthm}
\usepackage{appendix}
\usepackage{epstopdf}
\usepackage{booktabs,multirow}
\usepackage{float}
\usepackage{dsfont}
\usepackage{stmaryrd}
\usepackage{xcolor}
\usepackage{cases}

\newtheorem{remark}{Remark}
\newtheorem{example}{Example}

\journal{??}

\theoremstyle{plain}
\newtheorem{thm}{Theorem}[section]

\theoremstyle{remark}

\def\d{\textrm{d}}
\def\diag{\mathrm{diag}}

\def\bx{\mathbf{x}}

\def\bI{\mathbf{I}}

\def\bR{\mathbf{R}}

\def\bM{\mathbf{M}}

\def\cL{\mathcal{L}}

\def\cA{\mathcal{A}}
\def\cB{\mathcal{B}}
\def\cN{\mathcal{N}}

\newcommand{\ben}{\begin{eqnarray}}
\newcommand{\een}{\end{eqnarray}}
\newcommand{\beq}{\begin{equation}}
\newcommand{\eeq}{\end{equation}}
\newcommand{\bea}{\begin{array}}
\newcommand{\eea}{\end{array}}
\newcommand{\bef}{\begin{figure}[H]}
\newcommand{\eef}{\end{figure}}

\numberwithin{equation}{section}

\begin{document}

\begin{frontmatter}

\title{Supplementary Variable Method for Developing Structure-Preserving Numerical Approximations to Thermodynamically Consistent Partial Differential Equations}

\author[mymainaddress,mymainaddress2]{Yuezheng Gong}
\author[mysecondaryaddress,mymainaddress2]{Qi Hong}
\author[mythirdaryaddress]{Qi Wang}
\cortext[mycorrespondingauthor]{Corresponding author}
\ead{qwang@math.sc.edu}
\address[mymainaddress]{College of Science, Nanjing University of Aeronautics and Astronautics, Nanjing 210016, China}
\address[mysecondaryaddress]{Beijing Computational Science Research Center, Beijing 100193, China}
\address[mymainaddress2]{{Jiangsu Key Laboratory for Numerical Simulation of Large Scale Complex Systems}}
\address[mythirdaryaddress]{Department of Mathematics, University of South Carolina, Columbia, SC 29208, USA}
\begin{abstract}

We present a new temporal discretization paradigm  for developing energy-production-rate preserving numerical approximations to thermodynamically consistent partial differential equation systems, called the supplementary variable method.   The central idea behind it is to introduce a supplementary variable to the thermodynamically consistent model to make the over-determined equation system, consisting of the thermodynamically consistent PDE system, the energy definition and the energy dissipation equation, structurally stable. The supplementary variable allows one to retain the consistency between the energy dissipation equation and the PDE system after the temporal discretization. We illustrate the method using a dissipative gradient flow model. Among virtually infinite many possibilities, we present two ways to add the supplementary variable in the gradient flow model to develop energy-dissipation-rate preserving algorithms. Spatial discretizations are carried out using the pseudo-spectral method. We then compare the two new schemes with the energy stable SAV scheme and the fully implicit Crank-Nicolson scheme. The results favor the new schemes in the overall performance. This new numerical paradigm can be applied to any  thermodynamically consistent models.
\end{abstract}


\end{frontmatter}

\section{Introduction}
Nonequilibrium phenomena require dynamical models derived from laws and principles of nonequilibrium thermodynamics to describe. The laws of thermodynamics, especially, the second law of thermodynamics or the equivalent generalized Onsager principle are fundamental laws/principles for developing such models for nonequilibrium phenomena not far from equilibria \cite{Onsager1931a,Onsager1931b,Yang&Li&Forest&Wang2016}. Assuming the thermodynamic variables describing nonequilibrium phenomena of a certain system are denoted as  $\Phi\in {\bf R}^n$, where $n$ is the number of the variables, and the free energy of the system is described by a functional $F[\Phi]$ in isothermal cases (for nonisothermal cases, we use an entropy functional instead). The Onsager linear response theory yields the dynamical equation
\ben
\bR \cdot \dot{\Phi}=-\bM \cdot \mu,\quad \mu = \frac{\delta F}{\delta \Phi},\label{gfeq}
\een
where $\bM$ and $\bR$ are operators, $\frac{\delta F}{\delta \Phi}$ is the variation of $F$, $\bM^{-1} \bR$ is known as the friction operator and  $\bR^{-1} \bM$ the mobility operator when they exist. This equation system provides relaxation dynamics for the nonequilibirum  state to return to equilibria in dissipative systems or to oscillate in nondissipative systems. For simplicity, we consider the case of  $\bR=\bI$ in this paper so that $\bM$ is the mobility operator.

In general, the mobility operator can be decomposed into two parts:
\ben
\bM=\bM_{a}+\bM_{s},
\een
where $\bM_a$ is the antisymmetric part and $\bM_s$ is the symmetric part. The time rate of change of the free energy  is given by
\ben
\dfrac{\d F}{\d t}= -\int_{\Omega} \frac{\delta F}{\delta \Phi} \bM \frac{\delta F}{\delta \Phi} \textrm{d}\bx=-\int_{\Omega} \frac{\delta F}{\delta \Phi} \bM_s \frac{\delta F}{\delta \Phi} \textrm{d}\bx ,\label{energy-dissipation-rate}
\een
under proper boundary conditions on $\Phi$.
If $\bM_a=0$ and $\bM_s\geq 0$,  the system is called a dissipative system since $\d F/ \d t \leq 0$. While only $\bM_s=0$,  the system is a conservative system because $\d F/ \d t=0$.
Apparently, \eqref{energy-dissipation-rate} is an equation  deduced from \eqref{gfeq}. It can therefore be viewed as a consistent constraint for \eqref{gfeq}.  When approximating equation \eqref{gfeq} numerically, one would like to preserve energy dissipation property at the discrete level, i.e., to arrive at an approximate equation to \eqref{energy-dissipation-rate} that is derivable from the approximate equation of \eqref{gfeq}. A numerical algorithm of such a property is called an energy stable algorithm or scheme. 

In many thermodynamical models given in the form of \eqref{gfeq}, the total free energy is given by 
\beq
F = \frac{1}{2}(\Phi, \cL \Phi ) + \big( f (\Phi, \nabla \Phi), 1\big),\label{free-energy}
\eeq
where $(\bullet,\bullet)$ is the  inner product in a $L^2$ space,  $\cL$ is a linear, self-adjoint, positive definite operator and $f$ or $(f,1)$ is bounded below for the physically accessible states of $\Phi$. This is also known as the gradient flow model in the literature. In this  model, if we view $F$ as one of the thermodynamical variables, \eqref{gfeq}, \eqref{energy-dissipation-rate}, and \eqref{free-energy} constitute an {\it over-determined system} of equations with $n+2$ equations and $n+1$ unknowns $(\Phi, F)$. This system of equations is consistent and solvable should gradient flow model \eqref{gfeq} is solvable with free energy \eqref{free-energy} and proper initial and boundary conditions. However, the structural consistency among the equations in the system can be easily broken under small perturbations to the system, leading to a system that is not solvable because of inconsistency among the equations. We label this scenario in the over-determined system of equations as structural unstable.

Traditionally, to retain consistency among the equations in the over-determined system in numerical approximations, one  discretizes \eqref{gfeq} and \eqref{free-energy}, and then try to put together a consistent discretized energy dissipation equation or inequality for \eqref{energy-dissipation-rate} to arrive at energy stable algorithms. However, there is no guarantee that this approach would be  successful. Then, there came about the convex splitting  method for the gradient flow with a free energy given as a difference of two convex functions, in which {\it  one implements a mismatched temporal discretization on the two convex parts of the free energy, respectively, to achieve energy stability} \cite{Eyre1998}; and the stabilizer approach \cite{Yang09}, in which {\it one modifies the gradient flow equation by adding higher order perturbations}. All these approaches try to put together a consistent, discrete  energy dissipation equation  following the traditional approach mentioned above within the over-determined partial differential equation system.

Recently, the energy quadratization (EQ) approach coined by Xiaofeng Yang, Jia Zhao and Qi Wang following the work of Ganzales and Tierra on liquid crystal models, provides a systematic approach to derive energy-dissipation-rate preserving schemes for thermodynamically consistent models \cite{Guill2013On,Yang2017Numerical,Zhao2017Numerical}. In this approach, {\it one introduces new unknowns (known as auxiliary variables) to reformulate the energy into a quadratic form and in the meantime supplement the system with the same number of dynamical equations for the new auxiliary variables}. This process is known as the EQ reformulation, which effectively embeds the gradient flow model into a higher dimensional phase space with a quadratic free energy while retaining the thermodynamically consistent structure and property (including the energy dissipation property.) The EQ reformulated gradient flow model together with the definition of the quadratic free energy and the energy dissipation equation constitute {\it a reformulated  over-determined, consistent and solvable gradient flow system.} One designs linear, energy stable schemes from the EQ reformulated gradient flow model. Analogously, the SAV method follows the same idea \cite{shen2018scalar}.

We illustrate the idea using an example here.  We assume $f$ only depends on $\Phi$, but not its spatial derivatives. In this case, one introduces a $q$ variable to transform the free energy into a quadratic form:
\ben
q = \sqrt{2\left(f +\frac{A}{|\Omega|}\right)}, \quad
F = \frac{1}{2}(\Phi, \cL \Phi) + \frac{1}{2}\|q\|^2 - A,\label{q-f-eq}
\een
where  $A$ is a constant large enough to make $q$ well-defined.
Then, one takes the time derivative of the algebraic equation of $q$ to arrive at an equivalent dynamical equation of $q$. The EQ-reformulated  system in a higher dimensional phase space recasts into
\beq \label{eq:gradient-flow-EQ-scalar}
\left\{
\bea{l}
\frac{\partial}{\partial t}\Phi = -\bM \Big(\cL \Phi + q \frac{\partial q}{\partial \Phi} \Big), \\
\frac{\partial}{\partial t} q = \frac{\partial q}{\partial \Phi} \cdot \frac{\partial \Phi}{\partial t}.
\eea
\right.
\eeq
Denoting $\Psi = \left( \Phi,\; q \right)^T$, the above system is rewritten into \cite{Gong&Zhao&WangSISC2020}
\beq \label{eq:gradient-flow-EQ-full}
\frac{\partial}{\partial t} \Psi = -\cN \cB \Psi,
\eeq
where $\cN = \cA^* \bM \cA,$ $\cA^*$ is the adjoint operator of $\cA,$ $$\cA = \left (\bI_n \quad \frac{\partial q}{\partial \Phi} \right)_{n,n+1},\quad \cB = \diag(\cL, 1)_{n+1,n+1},$$  $\cB$ is a linear, self-adjoint, positive definite operator.
The energy dissipation rate is given by
\beq
\frac{\d F}{\d t} = \Big( \frac{\delta F}{\delta \Psi} , \frac{\partial\Psi}{\partial t} \Big) = -\Big( \cB \Psi , \cN  \cB\Psi \Big) \leq 0.\label{eq:energy-diss}
\eeq
Recently, Gong, Zhao and Wang showed that one can devise arbitrarily high order energy stable schemes to solve the over-determined system consisting of \eqref{q-f-eq}, \eqref{eq:gradient-flow-EQ-full} and \eqref{eq:energy-diss} \cite{Gong&Zhao&WangSISC2020,Gong&ZhaoAML2019,Gong&Zhao&WangCPC2020,Gong&Zhao&WangArXiv2019}. In contrast to the traditional method, the EQ approach guarantees the consistency and solvability in the discretized system due to the reformulated free energy is quadratic  whereas the traditional approach may not.
Most recently, a class of new methods rooted in the scalar auxiliary variable approach have been devised to ensure linearity as well as energy stability in the numerical schemes developed for gradient flow models \cite{AkrivisLiLi:2019,Yang&g2020}. 

Here, we take look at the issue of developing energy stable numerical approximations from a new perspective. We begin with the over-determined, consistent system of equations consisting of \eqref{gfeq}, \eqref{energy-dissipation-rate}, and \eqref{free-energy}. Now that this system is structurally unstable, it can easily lose the consistency among the equations subject to any small perturbations. How can we make it structurally stable. One way to make it structurally stable is to extend the system by adding supplementary variables to make it well-determined, i.e., the number of unknowns equals to the number of equations. There are virtually unlimited number of ways this can be done. To be consistent with the original problem, we accomplish this using perturbations and require that the well-determined, perturbed system include the original system as a special case. We would like to take a new point of view towards next present a  new paradigm to devise energy dissipation rate preserving schemes for the over-determined system consisting of \eqref{gfeq}, \eqref{energy-dissipation-rate} and \eqref{free-energy} to retain structural consistency and solvability in the over-determined system. We name it the supplementary variable method ({\bf SVM}). We will show that this includes the newly proposed Lagrange multiplier method  based on the SAV approach by Cheng, Liu and Shen as well as the generalized SAV method developed by Yang and Dong lately as special cases \cite{ChengA,Yang&g2020}. This method originates from the idea on how to enforce structural stability in the over-determined equation system by augmenting it with supplementary variables. Specifically, we perturb the system consisting of \eqref{gfeq}, \eqref{energy-dissipation-rate} and \eqref{free-energy} by adding a perturbation:
\ben
\left \{
\bea{l}
\dot{\Phi}=-\bM \cdot \mu,\quad
 \mu = \frac{\delta F}{\delta \Phi}+\alpha g[\Phi],\\
F = \frac{1}{2}(\Phi, \cL \Phi ) + \big( f (\Phi, \nabla \Phi), 1\big)+\alpha h(t),\\
\dfrac{\d F}{\d t}=-\int_{\Omega} \frac{\delta F}{\delta \Phi} \bM_s \frac{\delta F}{\delta \Phi}+\alpha j(t),
\eea\right.\label{extended}
\een
where $\alpha(t)(g[\Phi], h(t), j(t))^T$ is the perturbation, $\alpha(t)$ is a function of time,  $g[\Phi]$ is  a prescribed functional of $\Phi$, $h$ and $j$ are prescribed functions of t. $\alpha(t)$ is called the supplementary variable. If $\alpha(t)=0$, this perturbed  equation system reduces to the original, consistent and over-determined system. By adding $\alpha$, however,  we  lift the dynamical system into a higher dimensional phase space  to remove its over-determinedness and to effectively provide it with structural stability. The gained structural stability warrants us to obtain energy-dissipation-rate preserving numerical algorithms for \eqref{gfeq} under very general conditions so long as we keep the perturbation the same order as the truncation error when we discretize the system. We next detail the implementation of this method in one simple gradient flow case.

\section{Supplementary variable method--a new paradigm for developing energy-dissipation-rate preserving numerical algorithms}

 We explicitly rewrite system \eqref{gfeq}, \eqref{energy-dissipation-rate} as follows
\ben
&& \partial_t \Phi = -\bM \big( \cL\Phi + f'(\Phi) \big),\label{Phi_eq} \\
&& \mu = \cL\Phi + f'(\Phi),\label{mu_eq}\\
&& \frac{\textrm{d} F}{\textrm{d} t} = - ( \mu, \bM\mu).\label{F_eq}
\een
To derive energy-dissipation-rate preserving numerical approximations, we firstly  approximate the energy-dissipation-rate equation of the system given in  \eqref{Phi_eq}-\eqref{F_eq} as follows
\ben
&& \frac{\widetilde{\Phi}^{n+\frac{1}{2}}-\Phi^n}{\tau/2} = - \bM \left(\cL \widetilde{\Phi}^{n+\frac{1}{2}} + f'(\overline{\Phi}^{n+\frac{1}{2}})\right),  \label{Phi_tilde}\\
&& \mu^{*} = \cL \widetilde{\Phi}^{n+\frac{1}{2}} + f'(\widetilde{\Phi}^{n+\frac{1}{2}}),\label{mu_star}\\
&& \frac{\widetilde F^{n+1}-F[\Phi^n]}{\tau} = - ( \mu^*, \bM\mu^*),\label{energy-constraint}
\een
where $\tau$ is the time step,
\ben
\overline{\Phi}^{n+\frac{1}{2}} = \frac{3\Phi^n - \Phi^{n-1}}{2}, \quad F[\Phi^n] = \frac{1}{2}(\Phi^n, \cL \Phi^n ) + \big( f (\Phi^n), 1\big).
\een
It follows from \eqref{Phi_tilde} that
\ben
\label{Phi_tilde_ex} \widetilde{\Phi}^{n+\frac{1}{2}} = \left(1+\frac{\tau}{2}\bM\cL\right)^{-1} \Big( \Phi^n - \frac{\tau}{2}\bM f'(\overline{\Phi}^{n+\frac{1}{2}}) \Big), \quad
{\widetilde F^{n+1}}=F[\Phi^n]-{\tau} ( \mu^*, \bM\mu^*).
\een
This is a 2th order approximation to the energy-dissipation-rate equation with a truncation error in the order of $O(\tau^2)$.

Secondly, 
we modify \eqref{Phi_eq} by a time-dependent supplementary variable $\alpha(t)$ together with a user supplied functional $g[\Phi]$:
\beq \label{LagrangePDEs}
\left\{
\bea{l}
\partial_t \Phi = -\bM \big( \cL\Phi + f'(\Phi) \big) + \alpha g[\Phi], ~ t\in (t^n, t^{n+1}], \\
F[\Phi(t^{n+1})] = \widetilde F^{n+1},
\eea
\right.
\eeq
where $g[\Phi]$ is a given functional that may depend on $\Phi$ and its derivatives. There is a grate deal of flexibility in determining how to modify the gradient flow model with a supplementary variable. It is clearly an open problem for this approach.   We simply present two special implementations in this study to illustrate the idea.
\begin{itemize}
\item One is
	\beq\label{gc1} g[\Phi] =  \bM f'(\Phi),\eeq which leads to
	\beq\label{modified_eq1} \partial_t \Phi = -\bM \Big( \cL\Phi + (1-\alpha)f'(\Phi) \Big). \eeq
 Here, the chemical potential is modified into
 \beq
 \label{shen_mu} \mu =  \cL\Phi + (1-\alpha)f'(\Phi).
\eeq
\item
	The other is
	\beq
\label{gc2} g[\Phi] =   -\bM \Big(  \cL\Phi + f'(\Phi)  \Big),
\eeq which implies
	 \beq
\label{modified_eq2} \partial_t \Phi = -( 1 +\alpha ) \bM \big( \cL\Phi + f'(\Phi) \big).
\eeq
Here, the mobility  is modified into $(1+\alpha)\bM$.
\end{itemize}

We note that a system consisting of  \eqref{F_eq}, \eqref{modified_eq1} and \eqref{shen_mu}  is equivalent to the reformulated model proposed in \cite{ChengA} using what they call the Lagrange multiplier method. The supplementary variable $\alpha$ can also be viewed as a perturbation variable since when $\alpha(t)=0$, the modified/perturbed system reduces to the original one. {\it By supplementing the over-determined system with a new variable without introducing an equation, we effectively remove the over-determinedness by making the modified system a well-determined in a higher dimensional phase space! In addition, the modified model reduces to the original one at $\alpha=0$. } We reiterate that there is  a great deal of flexibility to place the supplementary variable in the gradient flow model, which differentiates it from the Lagrange multiplier method and the generalized SAV method.

Thirdly, we apply the implicit-explicit Crank-Nicolson scheme in time to \eqref{LagrangePDEs} to arrive at the following semi-discrete system
\beq \label{scheme}
\begin{cases}
	\delta_t^+ \Phi^n = - \bM \left(\cL \Phi^{n+\frac{1}{2}} + f'(\widetilde{\Phi}^{n+\frac{1}{2}})\right) + \alpha^{n+\frac{1}{2}} g[\widetilde{\Phi}^{n+\frac{1}{2}}],\\
	F[\Phi^{n+1}] = \widetilde F^{n+1}, \quad \delta_t^+ \Phi^n = \frac{\Phi^{n+1}-\Phi^n}{\tau}, \quad \Phi^{n+\frac{1}{2}} = \frac{\Phi^{n+1}+\Phi^n}{2},
\end{cases}
\eeq
where $\widetilde{\Phi}^{n+\frac{1}{2}},~ \widetilde F^{n+1}$ have been obtained from \eqref{Phi_tilde_ex}.
\begin{remark}
	We apply a special implicit-explicit Crank-Nicolson scheme to the system consisting of \eqref{modified_eq1}, \eqref{shen_mu} and \eqref{F_eq} and obtain the following discrete method
	 \beq \label{shen_scheme}
	 \begin{cases}
	 	\delta_t^+ \Phi^n = - \bM \mu^{n+\frac{1}{2}},\\
	 	\mu^{n+\frac{1}{2}} = \cL \Phi^{n+\frac{1}{2}} + (1-\alpha^{n+\frac{1}{2}}) f'(\overline{\Phi}^{n+\frac{1}{2}}),\\
	 	\dfrac{F[\Phi^{n+1}]-F[\Phi^n]}{\tau} = - ( \mu^{n+\frac{1}{2}}, \bM\mu^{n+\frac{1}{2}}).
	 \end{cases}
	 \eeq
	 It is not difficult to prove that the scheme \eqref{shen_scheme} is equivalent to the Lagrange multiplier approach proposed in \cite{ChengA}.
\end{remark}

Next we discuss how to solve system \eqref{scheme} efficiently . Let
\ben
&& \widehat{\Phi}^{n+1} = (1+\frac{\tau}{2}\bM\cL)^{-1} \Big( (1-\frac{\tau}{2}\bM\cL)\Phi^n - \tau \bM f'(\widetilde{\Phi}^{n+\frac{1}{2}}) \Big),\label{Phi_hat}\\
&& w^n = (1+\frac{\tau}{2}\bM\cL)^{-1} g[\widetilde{\Phi}^{n+\frac{1}{2}}],\label{w_eq}\\
&& \beta = \tau \alpha^{n+\frac{1}{2}}.
\een
It follows from the first equation of \eqref{scheme} that
\beq
\label{get_Phi} \Phi^{n+1} = \widehat{\Phi}^{n+1} + \beta w^n.
\eeq
Then we plug $\Phi^{n+1}$ into the second equation in \eqref{scheme} to obtain
\beq\label{beta_eq} F[\widehat{\Phi}^{n+1} + \beta w^n] = \widetilde F^{n+1}, \eeq
which is an algebraic equation for $\beta.$ In general, it can have multiple solutions, but one of them must be close to 0 and it approaches to zero as $\tau \to 0$. So, we solve for this solution using an iterative method such as the Newton iteration with 0 as the initial condition, it generally converges to a solution close to 0 when $\tau$ is not too large. After obtaining $\beta$, we update $\Phi^{n+1}$ using \eqref{get_Phi}. Following the work of Refs. \cite{calvo2006on,calvo2010projection}, the existence of solution $\beta$ is guaranteed under the conditions of the following theorem.

\begin{thm}\label{thm:order}
	If $\left( \frac{\delta F}{\delta \Phi}[\Phi^n] , g[\Phi^n] \right) \neq 0,$ there exists a $\tau^*>0$ such that \eqref{beta_eq} defines a unique function $\beta = \beta (\tau)$ for all $\tau\in[0,\tau^*]$ and scheme \eqref{scheme} with \eqref{Phi_tilde}-\eqref{energy-constraint} is of order $O(\tau^2)$.
\end{thm}

\begin{proof}
	For $\tau,\beta$ in a neighborhood of $(0,0)$, we define the real function
\ben
u(\tau,\beta) = F[\widehat{\Phi}^{n+1} + \beta w^n] -\widetilde  F^{n+1} = F[\widehat{\Phi}^{n+1} + \beta w^n] - F[\Phi^n] + \tau(\mu^* , \bM \mu^*),
\een
where $\mu^*$ is calculated in the first step. 
	It follows from  \eqref{mu_star}, \eqref{Phi_tilde_ex}, \eqref{Phi_hat} and \eqref{w_eq} that
\ben
u(0,0) = 0, \quad \frac{\partial u}{\partial \beta}(0,0) = \left( \frac{\delta F}{\delta \Phi}[\Phi^n] , g[\Phi^n] \right) \neq 0.
\een
	According to the implicit function theorem, there exists a $\tau^*>0$ such that  equation $u(\tau,\beta) = 0$ defines a unique smooth function $\beta = \beta (\tau)$ satisfying $\beta(0) = 0$ and $u\big(\tau,\beta(\tau)\big) = 0$ for all $\tau\in[0,\tau^*]$.
	
	Since $\widehat{\Phi}^{n+1}$ satisfies the following scheme
	\beq \frac{\widehat{\Phi}^{n+1}-\Phi^n}{\tau} = - \bM \left(\cL\frac{\widehat{\Phi}^{n+1}+\Phi^n}{2} + f'(\widetilde{\Phi}^{n+\frac{1}{2}})\right),\eeq
	through local truncation error analysis we have
	\ben
\widehat{\Phi}^{n+1} = \Phi(t^{n+1}) + O(\tau^3), \quad F[\widehat{\Phi}^{n+1} ] = F[\Phi(t^{n+1})] + O(\tau^3).
\een
In addition, we expand
\ben
u(\tau,\beta) = u(\tau,0) +  \beta \frac{\partial u}{\partial \beta}(\tau,0) + O(\beta^2),
\een
 with
	\ben
\bea{l}
u(\tau,0) = F[\widehat{\Phi}^{n+1}] - \widetilde F^{n+1} = O(\tau^3),\\[0.1cm]
\dfrac{\partial u}{\partial \beta}(\tau,0) = \dfrac{\partial u}{\partial \beta}(0,0) + O(\tau).
\eea
\een
	Then $\beta = \beta(\tau) = O(\tau^3)$ and the proposed scheme is of order $O(\tau^2)$.
\end{proof}

\begin{remark}
	If \eqref{gc2} is chosen, i.e. $g[\Phi] = -\bM\frac{\delta F}{\delta \Phi}[\Phi],$ then the sufficient condition in Theorem \ref{thm:order} reduces to $\left( \frac{\delta F}{\delta \Phi}[\Phi^n], \bM\frac{\delta F}{\delta \Phi}[\Phi^n] \right) \neq 0,$ which is usually satisfied when the steady state is not reached.
\end{remark}

\section{Numerical results}
We present a numerical example to demonstrate the practicability, accuracy, as well as energy stability of the proposed schemes. For simplicity, we use periodic boundary conditions for the numerical example below and name the resulting scheme {\bf SVM-I}  and {\bf SVM-II} when \eqref{gc1} and \eqref{gc2} are adopted, respectively.

\begin{example}[Cahn-Hilliard model]
	We consider the following Cahn-Hilliard phase field model with phase variable $\phi$
	\begin{align*}
	\partial_t \phi = \lambda \Delta \left(   - \varepsilon^2 \Delta \phi + \phi^3 - \phi  \right),
	\end{align*}
	where the free  energy is given by
	\begin{align*}
	F = \frac{ \varepsilon^2 }{2} \|  \nabla \phi \|^2 + \frac{1}{4}(  \phi^2 - 1 )^2.
	\end{align*}
	First of all, we present the  mesh refinement test in time to confirm the order of accuracy of the schemes. We set the computational domain as $\Omega = [0, 1]^2$ and the parameter values as $\varepsilon = 10^{-2}$ and $ \lambda = 10^{-3}$, respectively.
	This model is discretized spatially using a Fourier pseudo-spectral method with $256 \times 256$ spatial meshes.  We use initial condition $\phi = 0.25 \sin (2 \pi x) \cos(2 \pi y)$ and time steps $\tau  = 1.25e-2 \times \frac{1}{2^{k - 1}}$, $k = 1, 2, 3, \cdots$. The errors are calculated as the difference between the solution of the coarse time step and that of the adjacent finer time step.
	Both the discrete $L^2$ and $L^{\infty}$ errors of numerical solution $\phi$ at $t = 1$  are shown in Figure \ref{fig:test-order}, where  we observe that  the proposed schemes  yield second order convergence rates in time.
	\begin{figure}[H]
		\centering
		\subfigure[$L^2$ error.]{
			\includegraphics[width=0.35\textwidth,height=0.35\textwidth]{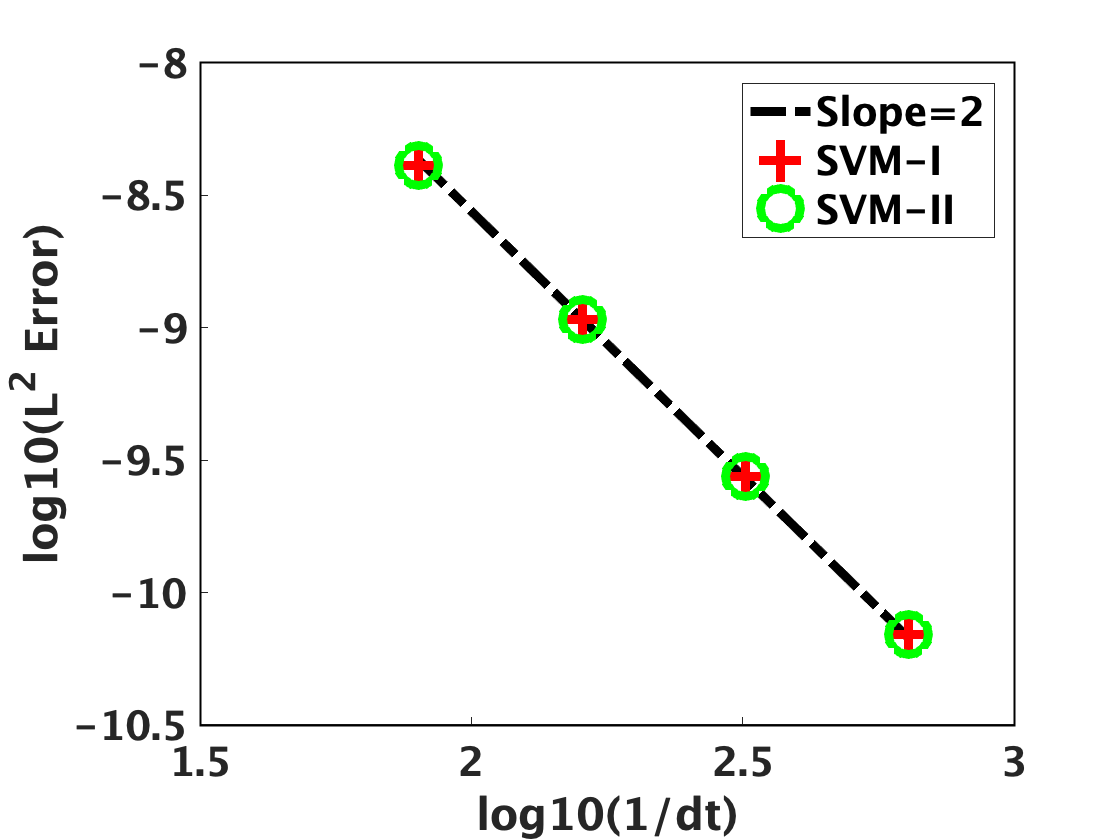}
		}\quad
		\subfigure[$L^{\infty}$ error.]{
			\includegraphics[width=0.35\textwidth,height=0.35\textwidth]{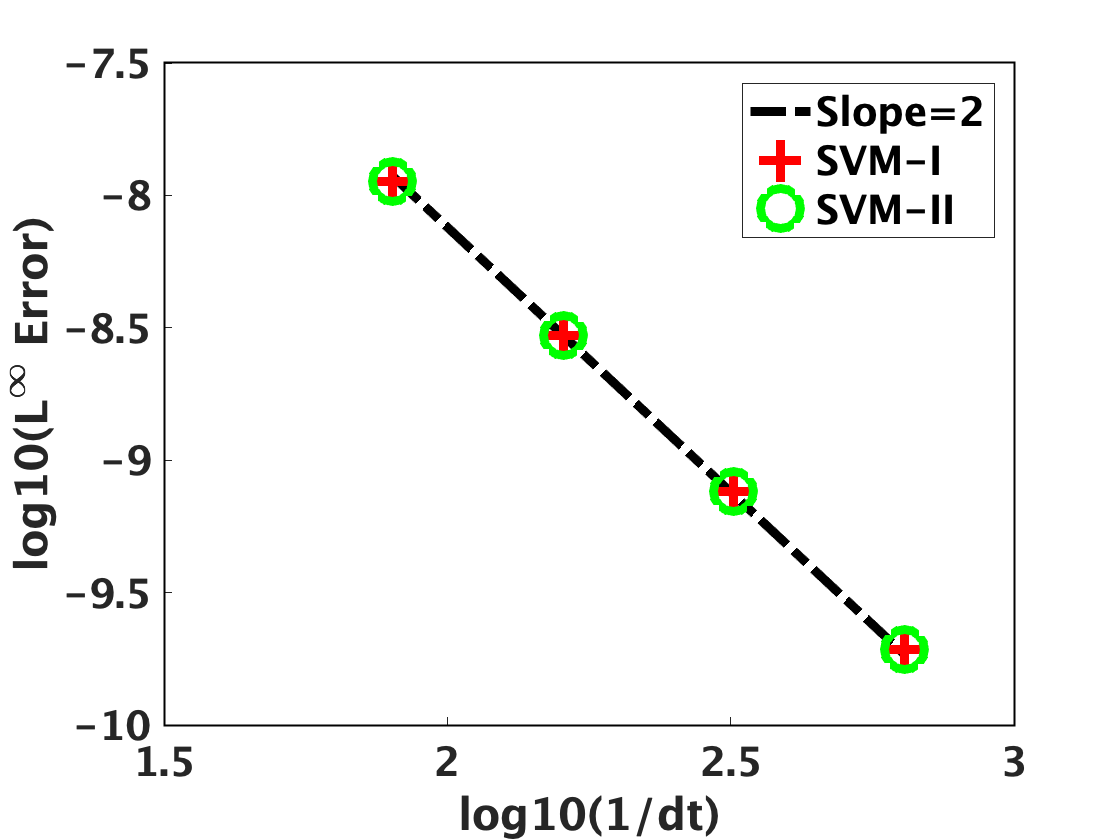}
		}
		\caption{Mesh refinement test in time for the two schemes. Here, we fix the number of spatial meshes at $N = 256$. Second order convergence rates are confirmed. } \label{fig:test-order}
	\end{figure}

	Next, we examine the computational efficiency of the two schemes: {\bf SVM-I} and {\bf SVM-II}.  We compare the two proposed schemes with the SAV scheme using the Crank-Nicolson method {\bf SAV-CN} \cite{shen2018scalar} and the fully implicit Crank-Nicolson scheme {\bf FICN} \cite{du1991numerical}. The result in the total CPU time to solve for  the solution using each scheme is summarized in Figure \ref{fig:test-CPU}. We observe that {\bf SVM-I} and {\bf SVM-II} for this model are less efficient than scheme {\bf SAV-CN}, but much more efficient than scheme {\bf FICN}. The price we pay using the proposed scheme in CPU time is that we have to solve a scalar nonlinear equation at each time step. In contrast, {\bf SAV} scheme solves a linear system while {\bf FICN} solves  a nonlinear system.
	
	\begin{figure}[H]
		\centering
		\subfigure{
			\includegraphics[width=0.60\textwidth,height=0.35\textwidth]{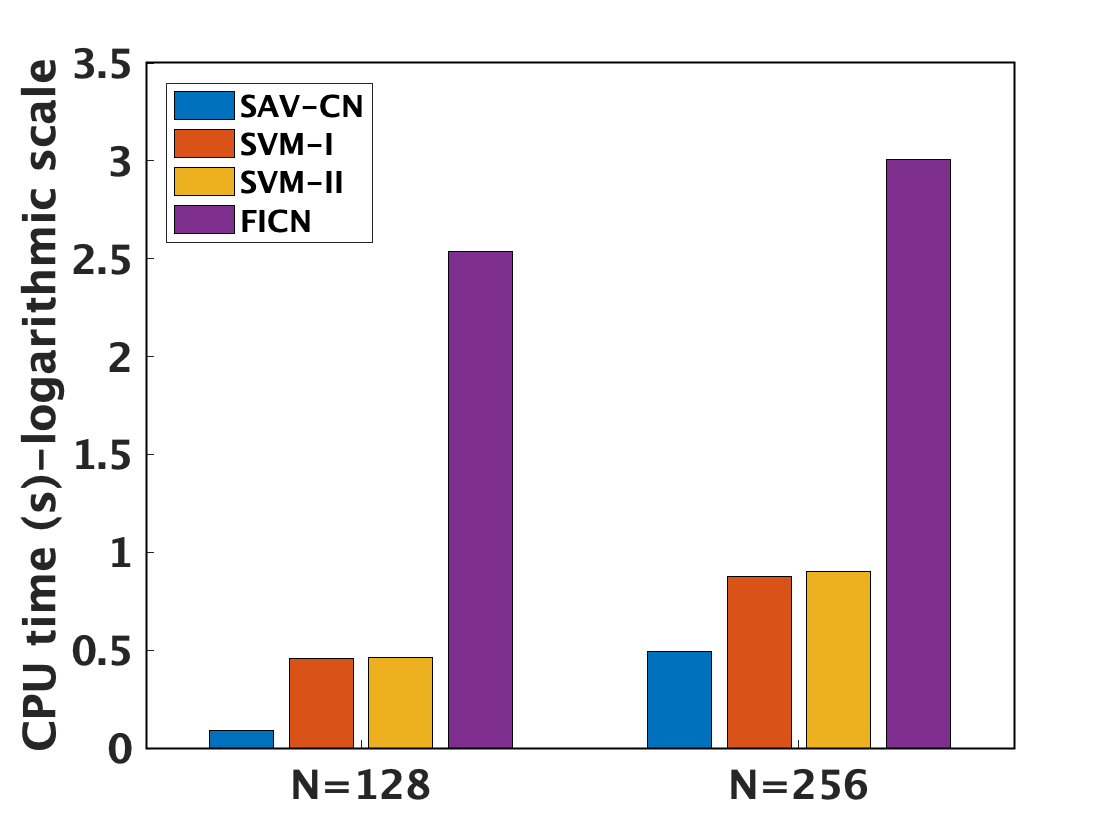}
		}
		\caption{Comparison of CPU time in logarithmic scale using the four different numerical methods with two sets of spatial meshes up to  $t = 1$, where the time step is chosen as $\tau = 1.0e-2$. The charts show that the two proposed schemes performs equally well, more slowly than {\bf SAV-CN}, but much faster than {\bf FICN}.} \label{fig:test-CPU}
	\end{figure}

Finally, we compare {\bf SVM-I}, {\bf SVM-II} and {\bf SAV-CN} in their accuracy in resolving the energy-dissipation-rate. We do it using a spatial discretization on $128^2$ meshes in $\Omega=[0,1]^2$. The parameter values are $\lambda = 1$ and $\varepsilon = 1.0e-2$.  We use the following initial condition \cite{Gong&Zhao&WangArXiv2019}
\begin{align*}
	\phi(x, y, 0) = 0.05 \left(  \cos(6\pi x) \cos(8 \pi y) + (\cos(8 \pi x) \cos(6 \pi y))^2 + \cos(2 \pi x -10 \pi y) \cos(4\pi x - 2\pi y)   \right).
\end{align*}
The initial profile undergoes a fast coarsening dynamics such that the algorithms must use small time steps in order to capture the correct coarsening dynamics.  The simulation results are depicted in Figure \ref{fig:test-fast-coarsening}, where the snapshots of $\phi(x,y ,t)$ at $t = 0.1$ are shown using different schemes and time steps.  We observe that scheme {\bf SAV-CN} predicts the correct solution at time step $\tau = 1.5625e-6$ but fails at time step $\tau = 3.125e-6$. For {\bf SVM-I}, it predicts correct numerical result  at time step  $\tau = 5.0e-5$, and {\bf SVM-II} performs even  better with time step $\tau = 2.0e-4$.
The errors in the total volume and the energy are plotted in Figure \ref{fig:test-energy-volume}. The numerical results show that both schemes conserve the total volume and capture the energy-dissipation-rate correctly using relatively larger time steps compared to the {\bf SAV} scheme. In addition, the supplementary variable $\alpha(t)$ remains close to zero except at a few initial time spots.

		\begin{figure}[H]
		\centering
				\subfigure[$\phi$ at $t = 0.1$ using {\bf SAV-CN} with time steps:  $\tau = 3.125e-6$, $\tau = 1.5625e-6$ and $\tau=1.0e-6$ (reference) (from left to right).]{
			\includegraphics[width=0.20\textwidth,height=0.20\textwidth]{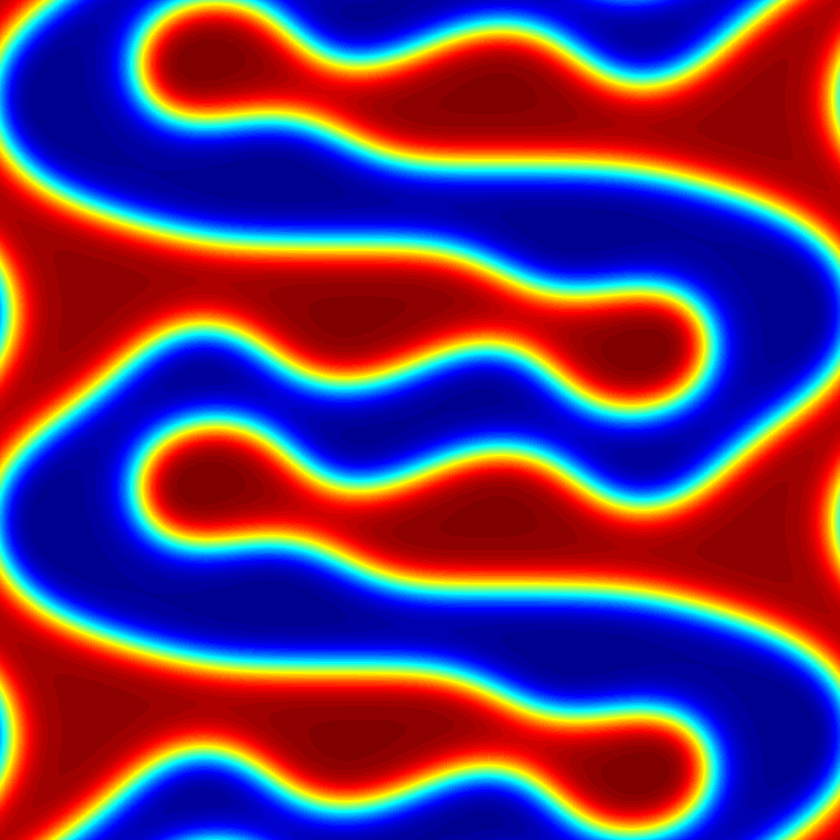}
			\includegraphics[width=0.20\textwidth,height=0.20\textwidth]{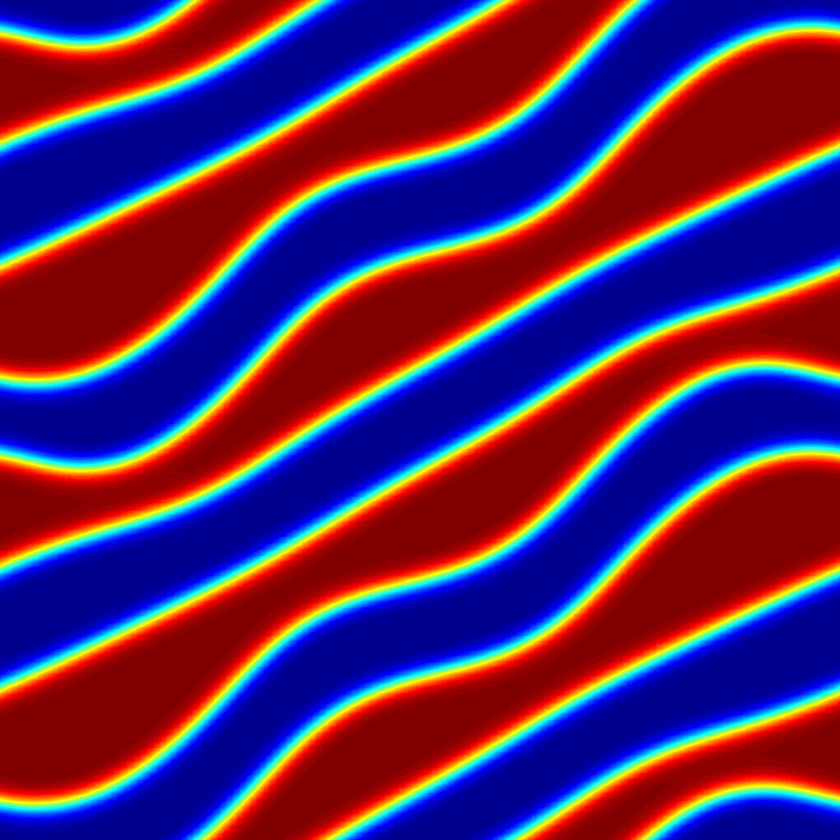}
			\includegraphics[width=0.20\textwidth,height=0.20\textwidth]{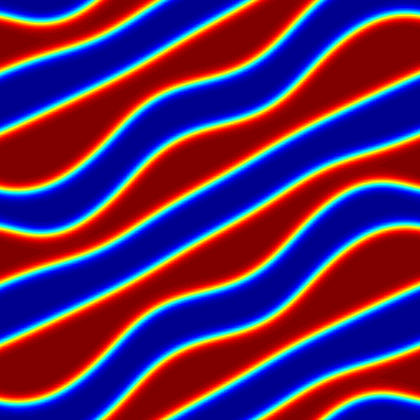}
		}
		\subfigure[ $\phi$ at $t = 0.1$ using {\bf SVM-I} with time steps:  $\tau = 5.0e-5$, $\tau = 4.0e-5$ and $\tau=1.0e-6$ (reference) (from left to right).]{
	       \includegraphics[width=0.20\textwidth,height=0.20\textwidth]{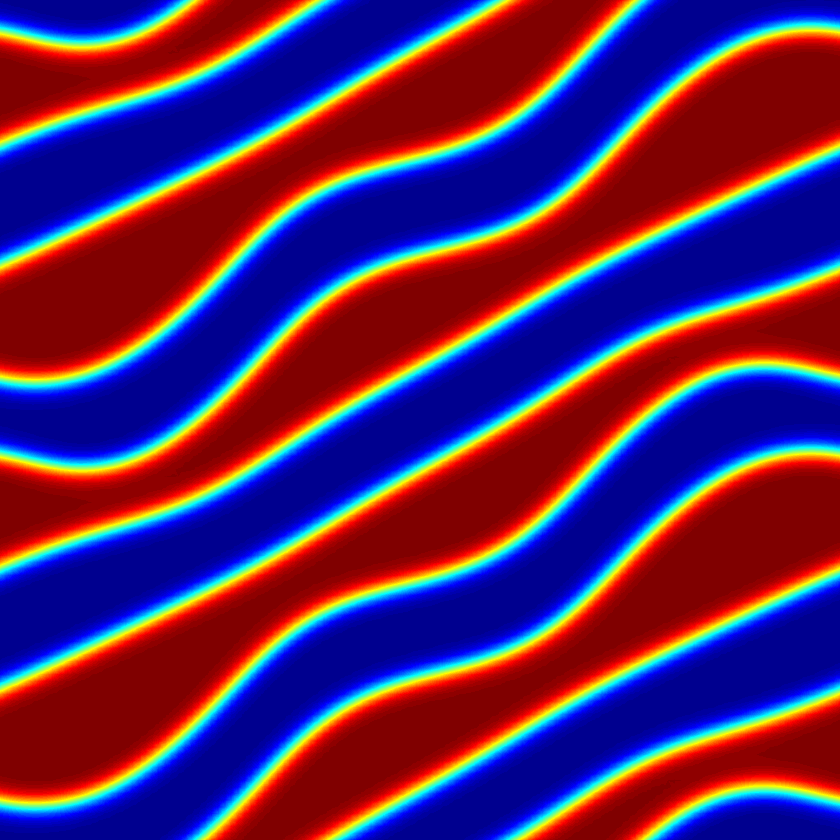}
			\includegraphics[width=0.20\textwidth,height=0.20\textwidth]{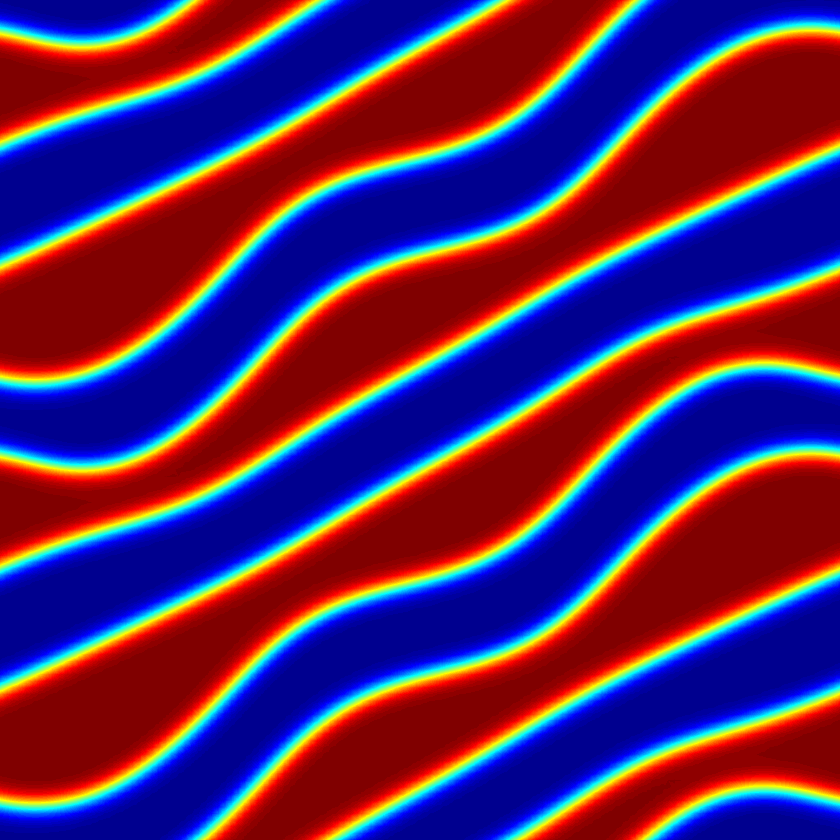}
			\includegraphics[width=0.20\textwidth,height=0.20\textwidth]{Figure/fast_coarsening/phi_reference.png}
		}
		\subfigure[$\phi$ at $t = 0.1$ using {\bf SVM-II} with time steps: $\tau = 2.0e-4$, $\tau = 1.5625e-4$ and $\tau=1.0e-6$ (reference) (from left to right).]{
	\includegraphics[width=0.20\textwidth,height=0.20\textwidth]{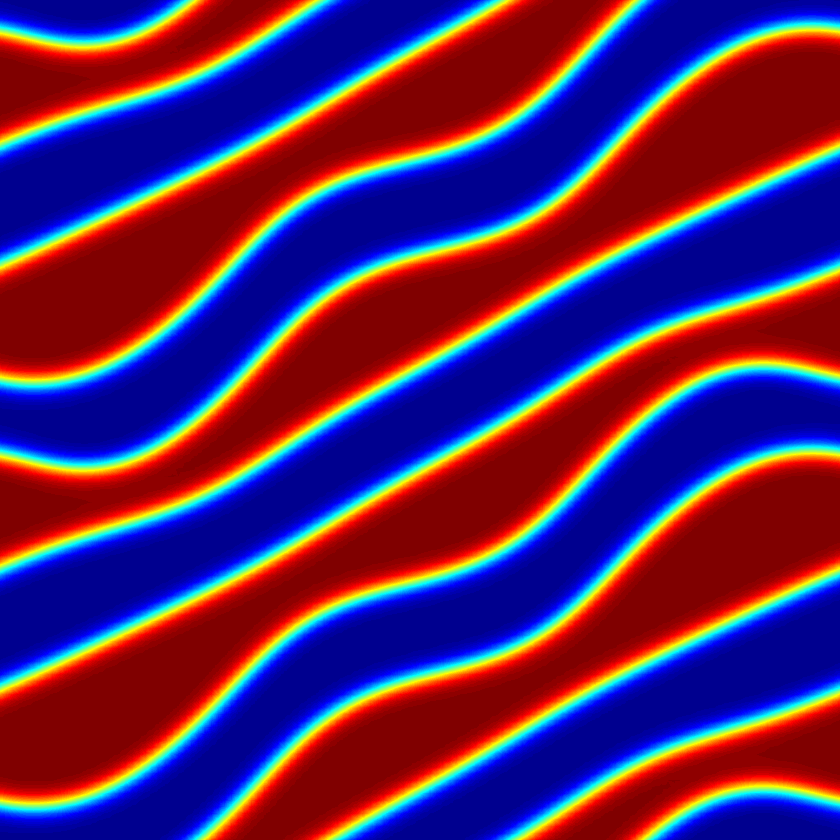}
	\includegraphics[width=0.20\textwidth,height=0.20\textwidth]{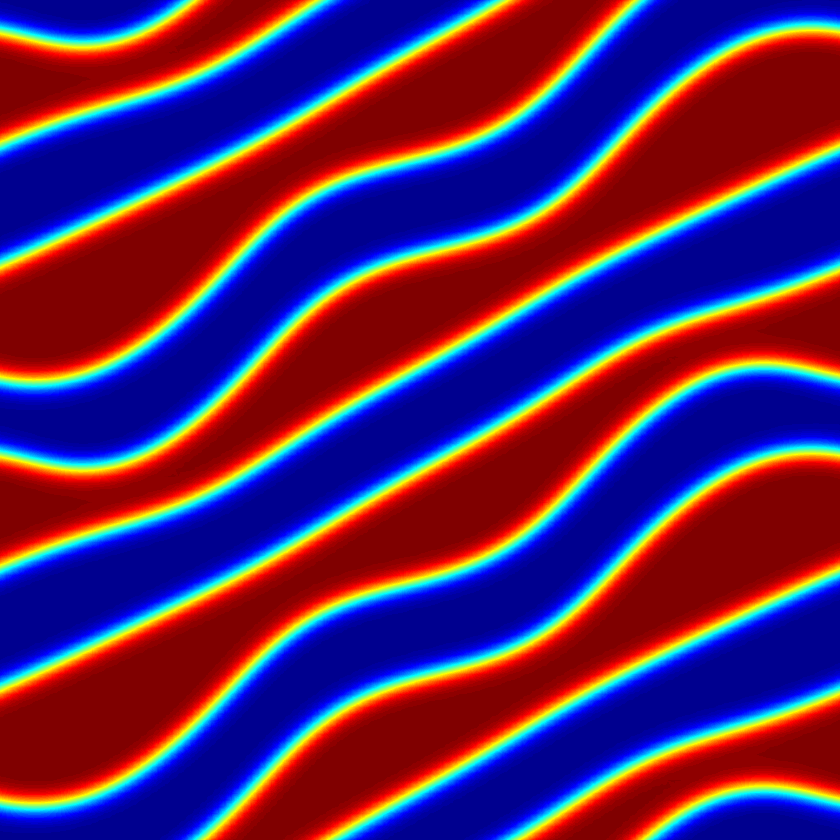}
	\includegraphics[width=0.20\textwidth,height=0.20\textwidth]{Figure/fast_coarsening/phi_reference.png}
}
		\caption{Comparison of three schemes at different time steps. The first sub-figure in each row indicates the ``maximum possible" time step to predict correct dynamics. These snapshots show that both schemes  {\bf SVM-I} and {\bf SVM-II} perform better than {\bf SAV-CN} and {\bf SVM-II} performs the best.  \label{fig:test-fast-coarsening} }
	\end{figure}

	\begin{figure}[H]
	\centering
	\subfigure[]{
	\includegraphics[width=0.30\textwidth,height=0.30\textwidth]{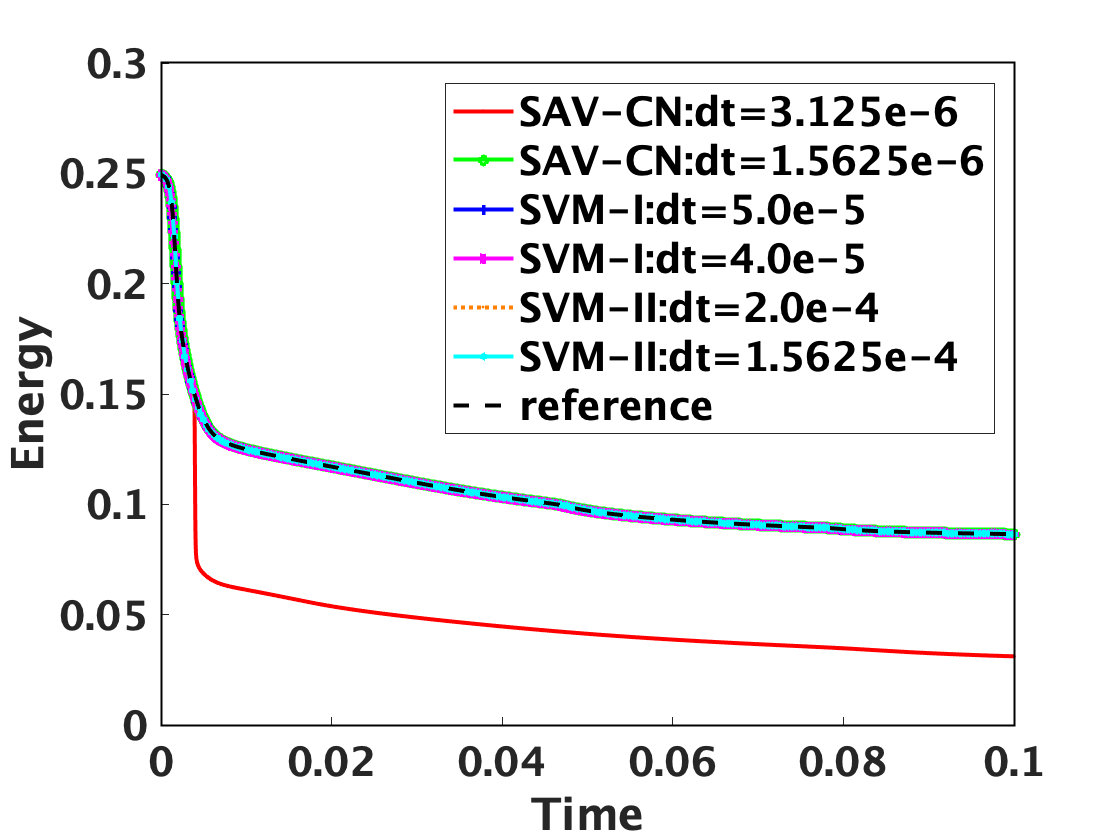}
}
	\subfigure[]{
		\includegraphics[width=0.30\textwidth,height=0.30\textwidth]{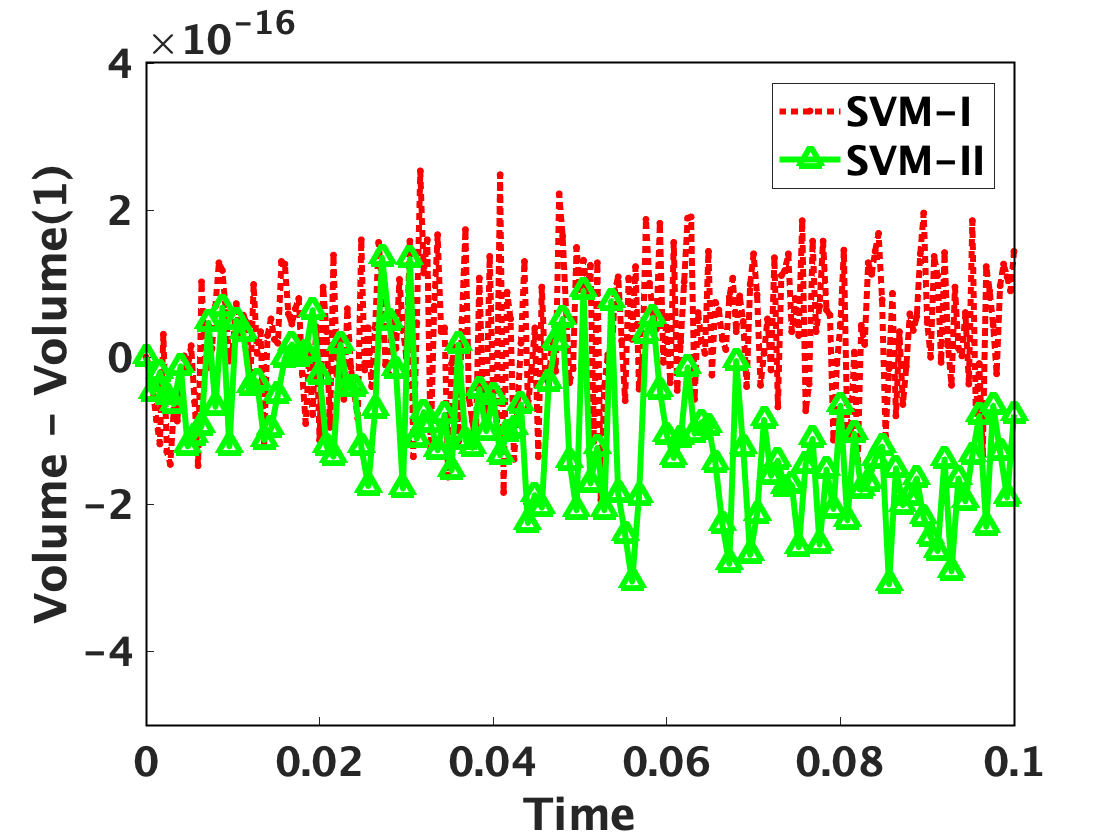}
	}
	\subfigure[]{
\includegraphics[width=0.30\textwidth,height=0.30\textwidth]{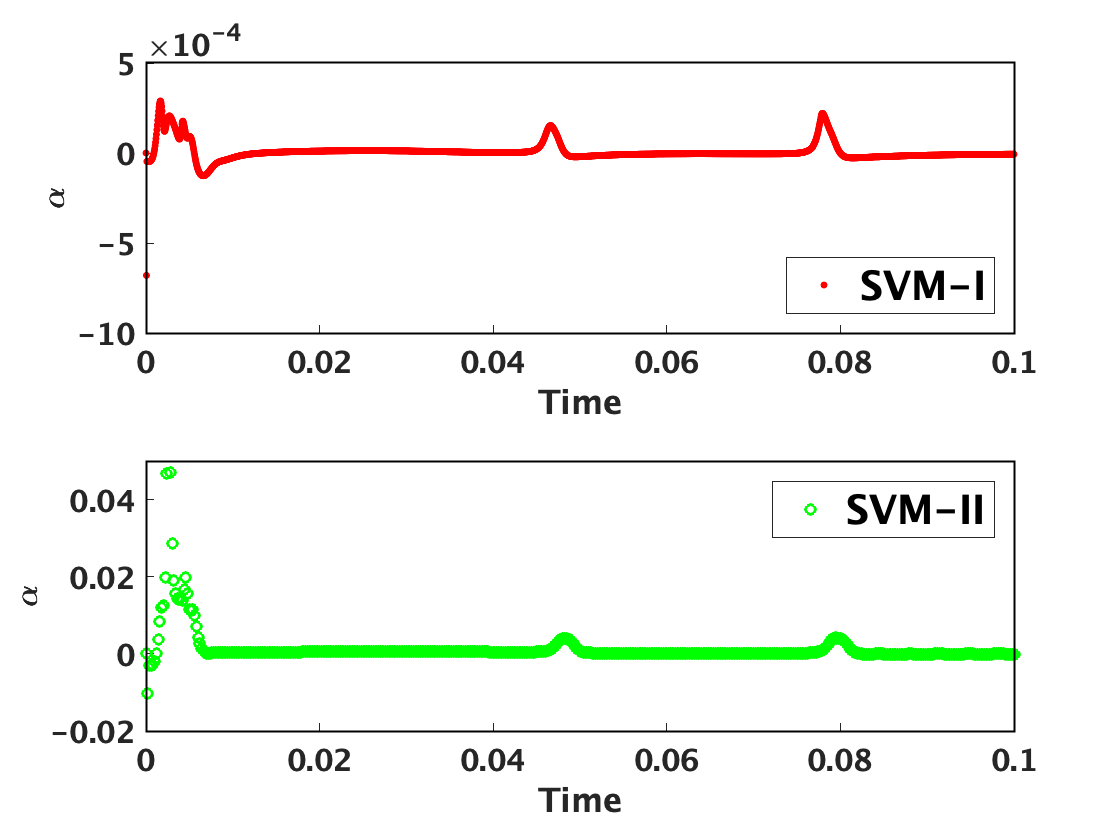}
}
	\caption{(a)Time evolution of the free energy computed using the three schemes with various time steps. The subfigure shows that {\bf SVM-I} and {\bf SVM-II} can work well with much larger time steps than {\bf SAC-CN} does, while {\bf SVM-II} performs the best. (b)
		Time evolution of the error in the total volume using {\bf SVM-I} and {\bf SVM-II} with $\tau = 5.0e-5$ and $\tau = 2.0e-4$, respectively.  The results show that the proposed schemes preserve the total volume very well. (c) The evolution of supplementary variable $\alpha$. This subfigure indicates  $\alpha$  may fluctuate near zero initially, but eventually, settles down close to zero.
		\label{fig:test-energy-volume}}
\end{figure}	
	\end{example}

\section{Conclusion}	
 The numerical results demonstrate the proposed schemes based on the supplementary variable approach can  predict fast coarsening dynamics of the Cahn-Hilliard model accurately and outperform  {\bf SAV-CN}  in solution accuracy and {\bf FICN} in efficiency.	Here we simply present two convenient implementations of supplementary variables for developing energy-dissipation-rate preserving schemes. There can be many other ways guided by this paradigm to achieve energy stability, better computational efficiency and accuracy. We expect to see more property preserving numerical algorithms developed for thermodynamically consistent models guided by this paradigm in the future.

\section*{Acknowledgment}
Research is partially supported by the Foundation of Jiangsu Key Laboratory for Numerical Simulation of Large Scale Complex Systems (202001, 202002), the Natural Science Foundation of Jiangsu Province (award BK20180413),  the National Natural Science Foundation of China (award 11801269 and  NSAF-U1930402) and National Science Foundation of US (award DMS-1815921 and  OIA-1655740) and a GEAR award from SC EPSCoR/IDeA Program.


\end{document}